\documentclass[11pt]{amsart}
\usepackage[marginratio=1:1,height=660pt,width=480pt,tmargin=70pt]{geometry}
\usepackage[english]{babel}
\usepackage[utf8]{inputenc}
\usepackage{amsfonts}
\usepackage[color,matrix, all, 2cell]{xy}
\usepackage{amsmath, amssymb}
\usepackage{graphicx}
\usepackage[font=small,labelfont=bf]{caption}
\usepackage{epstopdf}
\usepackage[pdfpagelabels,hyperindex]{hyperref}
\usepackage{xcolor}
\usepackage{amsthm}
\usepackage{float}
\usepackage{listings}
\usepackage{marginnote}
\usepackage{diagbox}
\usepackage{longtable}
\usepackage{mathrsfs}
\usepackage{lscape}
\usepackage{amsmath}
\usepackage{hyperref}
\usepackage{ dsfont }
\usepackage{mathtools}
\usepackage{float}
\usepackage{subfig}
\usepackage{hyperref}
\usepackage{amssymb}
\usepackage{nccmath}
\usepackage{empheq}
\usepackage{booktabs}
\usepackage{multirow}
\usepackage{siunitx}
\usepackage{gensymb}
\usepackage{upgreek}
\usepackage{mathtools}
\usepackage{derivative}
\usepackage{amsmath,amssymb}
\usepackage{witharrows}

\newcommand{\dd}{\mathrm{d}}
\newcommand{\longsquiggly}{\xymatrix{{}\ar@{~>}[r]&{}}}

\usepackage{tikz-cd}

\usepackage{yhmath}

\DeclareFontFamily{U}{mathx}{}
\DeclareFontShape{U}{mathx}{m}{n}{ <-> mathx10 }{}
\DeclareSymbolFont{mathx}{U}{mathx}{m}{n}
\DeclareFontSubstitution{U}{mathx}{m}{n}
\DeclareMathAccent{\widecheck}{0}{mathx}{"71}

\makeatletter
\newcommand{\xRrightarrow}[2][]{\ext@arrow 0359\Rrightarrowfill@{#1}{#2}}
\newcommand{\Rrightarrowfill@}{\arrowfill@\equiv\equiv\Rrightarrow}
\newcommand{\xLleftarrow}[2][]{\ext@arrow 3095\Lleftarrowfill@{#1}{#2}}
\newcommand{\Lleftarrowfill@}{\arrowfill@\Lleftarrow\equiv\equiv}
\makeatother

\newtheorem{thm}{Theorem}

\newtheorem{definition}{Definition}

\title{The CR geometry of the three-segment snake}
\author{Tymon Frelik}
\address{T.F.: Faculty of Physics\\ University of Warsaw\\Pasteura 5\\02-093 Warsaw\\ Poland and Center for Theoretical Physics\\ Polish Academy of Sciences\\ Al. Lotników 32/46\\ 02-668 Warsaw\\ Poland}
\email{t.frelik@student.uw.edu.pl}
\thanks{Support: The research leading to these results has received funding from the Norwegian Financial Mechanism 2014-2021 with project registration number 2019/34/H/ST1/00636. The author expresses gratitude for invaluable help and guidance of Paweł Nurowski and Katja Sagerschnig.}
\date{\today}
\begin{document}
\begin{abstract}
 We study the geometry associated with the kinematics of a planar robot known as the \emph{three-segment snake}, whose velocity distribution belongs to a class of $(2,3,5)$ distributions. We discover that, under certain assumptions on its construction parameters, the snake may be endowed with a CR structure of CR dimension 1 and real codimension 3. We solve the associated Cartan equivalence problem and find the invariants of the snake's CR structure.
\end{abstract}
\maketitle
\section{The three-segment snake}
\subsection{Introduction}
In the present note, we study the geometry associated with the kinematics of a planar robot we refer to as the \emph{three-segment snake} (also known as the \emph{three-link} or the \emph{three-edge snake}). This work is to be seen within a research program initiated by Paweł Nurowski, who observed that simple mechanical systems originating in robotics and control theory provide a rich reservoir for interesting geometric structures. The three-segment snake discussed here is the simplest such example.

Geometric treatment of the three-segment snake is found in the work of Masato Ishikawa \cite{Ish09}, who notably developed its control algorithm employing Ambrose-Singer’s holonomy theorem and verified it experimentally. Observing that the velocity distribution has the (small) growth vector $(2,3,5)$, Nurowski asked the following question: can the construction parameters of the snake, such as wheel positions and segment lengths, be adjusted in such a way that the algebra of infinitesimal symmetries is the 14-dimensional split real form of the complex exceptional Lie algebra $\mathfrak{g}_2$? 

The answer was obstructed by the computational complexity involved in finding the associated Cartan's quartic invariant \cite{Cartan1910}. In search of a geometric structure that could aid the \emph{hunt for a $G_2$-snake}, Nurowski observed \cite{Nur14} that, under a certain assumption on wheel placement, the $(2,3,5)$ distribution may be endowed with a natural CR structure, arising from the embedding of the 5-dimensional configuration space of the snake in the 8-dimensional space of all possible positions of all possible snakes in the Euclidean plane. 

The geometry of CR structures associated with embeddings of generic 5-manifolds in $\mathbb{C}^4$ has been studied by Merker et al. \cite{MPS14, MS15} and, independently, by Nurowski. In the present work, we study Nurowski's previously-unpublished results and give an invariant characterization of the CR structure associated with the snake by employing Cartan's equivalence method. 
\subsection{The configuration space}
The three-segment snake is a mobile robot moving in the Euclidean plane, constructed by attaching three line segments of fixed lengths. The line segments are free to rotate around the points of two connecting joints. Each segment has a pair of wheels attached at a fixed point. Their role is to direct the segments' motion.
\begin{figure}[]
    \centering
    \includegraphics[scale=0.56]{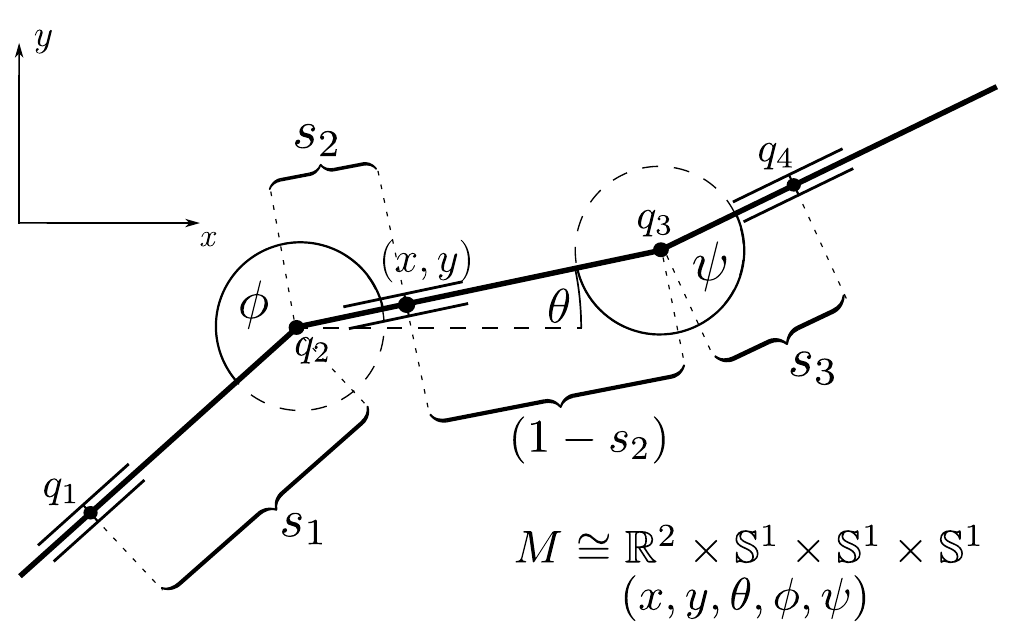}
    \caption{The three-segment snake in the Euclidean plane.}
    \label{fig:snakedrawing}
\end{figure}
The shape of a three-segment snake is described by four points $q_i=(x_i,y_i)\in\mathbb{R}^2,\;i=1,2,3,4$ (see Figure \ref{fig:snakedrawing}) so that any snake in any position may be described by a point in $q\in\mathbb{R}^8$. Without loss of generality, one may assume the length of the middle segment to be equal to 1 and the lengths of two outer segments to be determined by their wheels' placement. For any given snake, one requires the length of each segment to be constant, that is,
\begin{align*}
&h_1=|q_2-q_1|^2-s_1^2=(x_2-x_1)^2+(y_2-y_1)^2-s_1^2=0,\\
&h_2=|q_3-q_2|^2-1=(x_3-x_2)^2+(y_3-y_2)^2-1=0,\\
&h_3=|q_4-q_3|^2-s_3^2=(x_4-x_3)^2+(y_4-y_3)^2-s_3^2=0,
\end{align*}
where scaling parameters $s_1,s_3\in\mathbb{R}^+$ specify the length of two outer segments. A parameter $s_2\in ]0,1[$ describes the wheel placement on the middle segment at a point $s_2q_2+(1-s_2)q_3.$
These constraints are holonomic, since they depend only on positions. This exhibits the \emph{configuration space} of the snake, as a five-dimensional submanifold $M\subset\mathbb{R}^8$ given as the zero set
$$M=\{q\in\mathbb{R}^8 \;\;|\;\; h_1(q)=h_2(q)=h_3(q)=0\}.$$ 

The zero set of polynomials $h_1,h_2,h_3$ defines a foliation of $\mathbb{R}^8$ by five-dimensional leaves $M$, locally diffeomorphic to $\mathbb{R}^2\times \mathbb{S}^1\times \mathbb{S}^1\times \mathbb{S}^1$. 
Coordinates $(x,y,\theta,\phi,\psi)$ adapted to the foliation $M$ are most conveniently obtained by visual analysis of the system described in the Figure \ref{fig:snakedrawing}. This parametrization is manifestly invariant with respect to the action of the group $\mathrm{SE}(2):=\mathrm{SO(2)}\ltimes\mathbb{R}^2$ of orientation-preserving isometries of the plane. 
\section{Non-holonomic constraints on the snake's movement}
\label{nonholonomic}
The wheels are meant to prevent each of the segments from skidding sideways. This amounts to limiting the admissible displacements of each segment to a direction of the instantaneous orientation of the wheels. Formally, one requires that
$$ \dd q_{1}|| (q_{2}-q_{1}),\qquad\dd({(1-s_2)q_{2}+s_2q_{3}})|| (q_{3}-q_{2}),\qquad\dd q_{4}|| (q_{4}-q_{3}). $$
These conditions constitute a Pfaffian system, which may be rewritten in terms of vanishing of the forms $\Upsilon^1,\Upsilon^2,\Upsilon^3\in\Omega^1(\mathbb{R}^8)$ . In coordinates $(x,y,\theta,\phi,\psi)$, the system reads
\begin{align*}
&\Upsilon^1=\sin(\phi+\theta)\dd x-\cos(\phi+\theta)\dd y-(s_2\cos\phi-s_1)\dd\theta+s_1\dd\phi,\\& 
\Upsilon^2=\sin\theta\dd x - \cos\theta \dd y,\\
&\Upsilon^3=\sin(\psi+\theta)\dd x-\cos(\psi+\theta)\dd y -((1-s_2)\cos\psi-s_3)\dd\theta-s_3\dd\psi.
\end{align*}
The kernel of $\{\Upsilon^1,\Upsilon^2,\Upsilon^3\}$ is a rank 2 distribution $\mathcal{D}=\mathrm{span}_{C^\infty(M)}\{\xi_4, \xi_5\}\subset TM$ of admissible velocities with  
\begin{align*}
&\xi_4=   \partial_\theta-(1-\tfrac{s_2}{s_1}\cos\phi)\partial_\phi-(1-\tfrac{1-s_2}{s_3}\cos\psi)\partial_\psi,\\&\xi_5=\cos\theta\partial_x+\sin\theta\partial_y+\tfrac{1}{s_1}\sin\phi\partial_\phi-\tfrac{1}{s_3}\sin\psi\partial_\psi.
\end{align*}
The distribution $\mathcal{D}$ is bracket generating, since 
\begin{align*}&\xi_{3}:=[\xi_5,\xi_4]=\sin\theta\partial_x-\cos\theta\partial_y-\tfrac{1}{s_1}\left(\tfrac{s_2}{s_1}-\cos\phi\right)\partial_\phi+\tfrac{1}{s_3}\left(\tfrac{1-s_2}{s_3}-\cos\psi\right)\partial_\psi,
\\&\xi_{2}:=[\xi_5,\xi_{3}]=-\tfrac{1}{s_1^2}(1-\tfrac{s_2}{s_1}\cos\phi)\partial_\phi-\tfrac{1}{s_3^2}(1-\tfrac{1-s_2}{s_3}\cos\psi)\partial_\psi,\\
&\xi_{1}:=[\xi_4,\xi_{3}]=\cos\theta\partial_x+\sin\theta\partial_y+\tfrac{s_1^2-s_2^2}{s_1^3}\sin\phi\partial_\phi+\tfrac{(1-s_2)^2-s_3^2}{s_3^3}\sin\psi\partial_\psi
\end{align*}
satisfy $\wedge_{i=1}^5\xi_i\neq 0$. The tangent bundle is seen to admit a two-step filtration
$\mathcal{D}\subset[\mathcal{D},\mathcal{D}]\subset[\mathcal{D},[\mathcal{D},\mathcal{D}]]= TM,$
placing $\mathcal{D}$ in the class of distributions with the (small) growth vector $(2,3,5)$ or \emph{$(2,3,5)$ distributions}, for short. These are generic rank 2 distributions on 5-manifolds. This is in contrast with other possible nonholonomic planar systems, such as, for instance, a car \cite{HN19} pulling a single trailer \cite{Jea96}, whose growth vector is $(2,3,4,5)$.

We will call two $(2,3,5)$ distributions \emph{equivalent} if there exists a diffeomorphism $f:M\rightarrow M$ such that $f_*\mathcal{D}=\mathcal{D}$.

Generic rank 2 distributions in dimension 5 were first studies by Élie Cartan in the celebrated 1910 ``five-variables paper'' \cite{Cartan1910}. In the terminology of parabolic geometry \cite{CS09}, $(2,3,5)$ distributions are described by parabolic geometries of type $(\mathrm{G}_2,\mathrm{P})$ with $\mathrm{P}$ a 9-dimensional parabolic subgroup. In particular, the maximal symmetry algebra of a $(2,3,5)$ distribution is the split real form of the 14-dimensional exceptional Lie algebra $\mathfrak{g}_2$. This occurs for a flat (in the Cartan-theoretic sense) $(2,3,5)$ distribution.

In our case, only three infinitesimal symmetries $\varsigma_1,\varsigma_2,\varsigma_3$ of $\mathcal{D}$ are immediate:
$$\varsigma_1=\partial_x,\qquad\varsigma_2=\partial_y,\qquad\varsigma_3=y\partial_x-x\partial_y-\partial_\theta.$$
These are the generators of $\mathfrak{se}(2)=\mathfrak{so}(2)\oplus\mathbb{R}^2$.

In another modern treatment of $(2,3,5)$ distributions \cite{Nur05}, it is shown that each equivalence class of $(2,3,5)$ distributions defines a canonical conformal class of metrics of signature $(2,3)$, which can be algorhithmically computed. Furthermore, the Weyl tensor of such a conformal class contains the information about all basic differential invariants of $(M,\mathcal{D})$. The vanishing of the simplest of these invariants, the so-called \emph{Cartan's quartic}, is required for a maximally symmetric $(2,3,5)$ distribution. 

The question motivating the author's consideration of the three-segment snake robot model is the following. Does there exist a choice of parameters $s_1,s_2,s_3$ such that the algebra of symmetries of $\mathcal{D}$ is $\mathfrak{g}_2$? Even with the above algorithm, the question is too computationally complex to answer. We thus search for more geometric structure associated with the snake, which could reduce the computational complexity involved.
\section{The CR structure of $M\subset \mathbb{C}^4$}
Let $\mathcal{J}$ denote a complex structure on $\mathbb{R}^{2(n+k)}$ so that $(\mathbb{R}^{2(n+k)},\mathcal{J})\cong\mathbb{C}^{n+k}$. An embedding of a real $(2k+n)$-dimensional manifold $M$ in $\mathbb{C}^{n+k}$ gives rise to a CR structure of real codimension $n$ and CR dimension $k$ (type $(n,k)$). Such a CR structure may be described intrinsically by the distribution $\mathcal{D}:=\mathcal{J}(TM)\cap TM$ of rank $2k$. Then, $\mathcal{J}$ restricts to an endomorphism $\mathcal{J}:\mathcal{D}\rightarrow \mathcal{D}$. Such defined triple $(M,\mathcal{D},\mathcal{J})$ is referred to as a \emph{CR manifold}.

Recall that $M$ is naturally embedded in $\mathbb{R}^8$ as a submanifold of codimension 3. Furthermore, it is endowed with a corank 3 distribution $\mathcal{D}$. There are many ways of identifying $\mathbb{R}^8\cong\mathbb{C}^4$. In fact, the space of (linear) complex structures on $\mathbb{R}^8$ may be identified with $\mathrm{GL}(8,\mathbb{R})/\mathrm{GL}(4,\mathbb{C})$. Consequently, there are many ways that one could endow $M\subset\mathbb{R}^8$ with a CR structure, depending on a choice of a complex structure on $\mathbb{R}^8$. 

Hence, one may pose the following question. Does there exist a constant linear map $\mathcal{J}:\mathbb{R}^8\rightarrow \mathbb{R}^8$ satisfying $\mathcal{J}^2=-\mathrm{id}$ such that $\mathcal{J}(TM)\cap TM=\mathcal{D}$? The following theorem, originally proved by Paweł Nurowski \cite{Nur14}, provides a surprising answer.

\begin{thm}
    A complex structure $\mathcal{J}:\mathbb{R}^8\rightarrow\mathbb{R}^8$ such that $\mathcal{J}(TM)\cap TM=\mathcal{D}$ exists if and only if $s_2=\tfrac12$, i.e., the wheels of the middle segment are precisely in the center. Moreover, when $s_2=\tfrac12$, such $\mathcal{J}$ is unique.
\end{thm}
Therefore, a complex structure $\mathcal{J}$ endowing $(M,\mathcal{D},\mathcal{J})$ with a structure of a CR manifold exists if and only if $s_2=\tfrac12$. Then, the complex structure $\mathcal{J}$ reads
$$\mathcal{J}=\begin{psmallmatrix}
     0 & 0 & 0 & 0&|& -1 & 1 & 0 & 0 \\
    0 & 0 & 0 & 0& |& 0 & 1 & 0 & 0 \\
    0 & 0 & 0 & 0&  |& 0 & 0 & -1 & 0\\
    0 & 0 & 0 & 0&|& 0 & 0 & -1 & 1 \\
    -&-&-&-&-&-&-&-&-\\
     1 & -1 & 0 & 0&| & 0 & 0 & 0 & 0 \\
    0 & -1 & 0 & 0 &|& 0 & 0 & 0 & 0\\
     0 & 0 & 1 & 0&|&  0 & 0 & 0 & 0\\
     0 & 0 & 1 & -1 &|& 0 & 0 & 0 & 0
\end{psmallmatrix}.$$
The holomorphic coordinates $(z_1,z_2,z_3,z_4)$ on $\mathbb{C}^4$ defined by $\mathcal{J}$ are 
$$z_1=x_1+i(y_2-y_1), \quad z_2=x_2+iy_2, \quad z_3=x_3+iy_3,\quad z_4=x_4+i(y_4-y_3).$$

The complex structure $\mathcal{J}$ restricts to the real 2-plane field $\mathcal{D}$ and, at a point in $M$, is equivalent to the unique complex structure given by $i$ multiplication on $\mathbb{C}\cong\mathbb{R}^2$. The distribution $\mathcal{D}$ is spanned by $\zeta_\pm$, the $\pm i$-eigenvector fields of $\mathcal{J}$. 

\section{The Cartan equivalence problem for $(2,3,5)$ CR structures}
The CR structure $(M,\mathcal{D},\mathcal{J})$ associated with the snake will be considered among all generic type $(3,1)$ CR structures by employing Cartan's equivalence method. We will also refer to these structures as \emph{$(2,3,5)$ CR structures} for short. Their equivalence problem has been solved in \cite{MS15} and a canonical Cartan connection has been constructed in \cite{MPS14}. In this case, since one wishes to preserve the complex structure, we speak of biholomorphic rather than just diffeomorphic equivalence.

Generic CR structures of type $(3,1)$ may be realized as Cartan geometries of type $(G^7,\mathbb{C}^\times)$, where $G^7$ is a certain 7-dimensional real Lie group and $\mathbb{C}^\times=\mathrm{GL}(1,\mathbb{C})$ the multiplicative group. The corresponding 5-dimensional Lie algebra $\mathfrak{n}=\mathfrak{g}^7/\mathbb{C}$ is nilpotent and defined by the only three non-vanishing relations $$[e_1,e_2]=e_3,\quad [e_1,e_3]=e_4,\quad[e_2,e_3]=e_5$$ 
between its generators $e_1,e_2,e_3,e_4,e_5$.

\begin{definition}A coframe $\{\omega_1,\omega_2,\omega_3,\omega_4,\omega_5\}$ is called \emph{adapted to a $(2,3,5)$ CR structure $(\tilde{M},\tilde{\mathcal{D}},\tilde{\mathcal{J})}$} if $\tilde{\mathcal{D}}=\mathrm{ker}(\omega_1,\omega_2,\omega_3)$, $\omega_1=\bar\omega_2$ and $\omega_4=\bar\omega_5$, where conjugates are taken with respect to $\mathcal{J}$.\end{definition}

To restrict ourselves to five variables adapted to $M$ in solving the equivalence problem for $(M,\mathcal{D},\mathcal{J})$, we solve the system of polynomials $h_1,h_2,h_3$, substituting
$$2 z_1 - z_2 - \bar z_2 = 2 s_1 e^{i\beta_1}, \qquad z_2-z_3=e^{i\beta_2},\qquad 2 z_4 - z_3 - \bar z_3 = 2 s_3 e^{i\beta_3}$$
so that $(\beta_1,\beta_2,\beta_3,z_2,\bar{z}_{2})$ are the coordinates on $M$. To supplement $\Upsilon^1,\Upsilon^2,\Upsilon^3$ to an adapted coframe on $M$, it suffices to choose any pair of complex conjugate coordinate 1-forms, so we pick $\dd z_2, \dd\bar z_2$.
The coframe takes the form
$$\tilde\omega^1=\Upsilon^1+i\Upsilon^2,\qquad\tilde\omega^1=\Upsilon^1-i\Upsilon^2, \qquad \tilde\omega^3=\Upsilon^3,\qquad\tilde\omega^4=\dd z_2,\qquad\tilde\omega^5=\dd \bar{z}_2$$
and we verify that $\tilde\omega^1\wedge\tilde\omega^2\wedge\tilde\omega^3\wedge\tilde\omega^4\wedge\tilde\omega^5\neq0.$ 

The following two theorems are based on Nurowski's unpublished results on $(2,3,5)$ CR structures. Their equivalent alternative versions are in \cite{MPS14,MS15}.

\begin{thm}
Two $(2,3,5)$ CR structures are equivalent if their adapted coframes are related by a transformation belonging to the 10-dimensional $\mathrm{GL}(5,\mathbb{C})$-subgroup
$$H=\left\{\begin{psmallmatrix}
    g_1\bar{g_1}^2 & 0 & 0 & 0 & 0\\
    0 & g_1^2\bar{g_1} & 0 & 0 & 0\\
    \bar{g_3} & g_3 & g_1\bar{g_1} & 0 & 0\\
    \bar{g_5} & g_4 & \bar{g_2} & \bar{g_1} & 0\\
     \bar{g_4} & g_5 & g_2 & 0 & g_1
\end{psmallmatrix}\in \mathrm{GL}(5,\mathbb{C})\;:\; g_i\in \mathbb{C}, \;g_1\neq 0\right\},$$
given with respect to $\{\tilde{\omega}_1,\tilde{\omega}_2,\tilde{\omega}_3,\tilde{\omega}_4,\tilde{\omega}_5\}$.\end{thm}

\begin{thm}Given a coframe of 1-forms $\{\tilde{\omega}_1,\tilde{\omega}_2,\tilde{\omega}_3,\tilde{\omega}_4,\tilde{\omega}_5\}$ on $M$, one can always find $h\in H$, such that $\omega^i=h^i_j\tilde{\omega}^j$ satisfy the exterior differential system
\begin{align*}
\begin{split}
    & \dd \omega^1=-\omega^1\wedge\Omega^1+\mathbf{J}\omega^2\wedge\omega^4+\omega^3\wedge\omega^5,\\
    &\dd\omega^3=i\mathbf{T}\omega^1\wedge\omega^2+\mathbf{S}\omega^1\wedge\omega^3+\mathbf{L}\omega^1\wedge\omega^5+\mathbf{\bar S}\omega^2\wedge\omega^3\\&\qquad\;\;+\mathbf{\bar L}\omega^2\wedge\omega^4-\tfrac{1}{3}\omega^3\wedge\Omega^1-\tfrac{1}{3}\omega^3\wedge\Omega^2+i\omega^4\wedge\omega^5,\\
    &\dd\omega^4=\mathbf{Q}\omega^1\wedge\omega^2+\mathbf{G}\omega^1\wedge\omega^3-(\mathbf{S}-\mathbf{\bar J\bar V})\omega^1\wedge\omega^4-\mathbf{N}\omega^1\wedge\omega^5+\mathbf{K}\omega^2\wedge\omega^3-\mathbf{F}\omega^2\wedge\omega^4\\&\qquad\;\;+\mathbf{B}\omega^2\wedge\omega^5-\mathbf{V}\omega^3\wedge\omega^5- i\mathbf{A}\omega^3\wedge\omega^4+\tfrac{1}{3}\omega^4\wedge\Omega^1-\tfrac{2}{3}\omega^4\wedge\Omega^2,
    \end{split}
\end{align*}
where $\{\omega^1,\omega^2,\omega^3,\omega^4,\omega^5,\Omega^1,\Omega^2\}$, with $\omega^2=\bar\omega^1$, $\omega^5=\bar\omega^4$ and $\Omega^2=\bar\Omega^1$, constitute a coframe on $G^7$, while $\mathbf{J,T,S, L, Q, G,V,N,K,F,B,A}\in C^\infty(M)$. The latter 12 functions are the differential invariants of the CR structure and satisfy (among other similar relations not included here) the relevant relation \begin{align}\label{eq:firstinvariant}
    \dd \mathbf{J}=-(\mathbf{\bar N}-i\mathbf{ A J})\omega_3-\mathbf{\bar L}\omega_5+\tfrac{4}{3}\mathbf{J}\Omega_1-\tfrac{5}{3}\mathbf{J}\Omega_2+J_1\omega_1+J_2\omega_2+J_4\omega_4
\end{align}
with $J_1,J_2,J_3\in C^\infty(M)$. The (non-)vanishing of each of the 12 functions is an \emph{invariant property} of the corresponding CR structure.\end{thm}

For the CR structure $(M,\mathcal{D},\mathcal{J})$ of the snake, the primary invariant $\mathbf{J}$ (denoted by $R$ in \cite{MS15}) vanishes identically, which by \ref{eq:firstinvariant} also implies the vanishing of $\mathbf{N}$ and $\mathbf{L}$. As indicated in \cite{MS15}, $(2,3,5)$ CR structures belonging to the branch $\mathbf{J}\equiv 0$ admit a reduction of the $H$-structure to a 2-dimensional subgroup $H_{\mathbf{J}}\subset H$. In our case, it is determined to be of the form 
$$H_{\mathbf{J}}=\left\{\begin{psmallmatrix}
    g_1\bar{g_1}^2 & 0 & 0 & 0 & 0\\
    0 & g_1^2\bar{g_1} & 0 & 0 & 0\\
   0 & 0 & g_1\bar{g_1} & 0 & 0\\
    0 & 0 & g_2 & \bar{g_1} & 0\\
     0 & 0 & \bar{g_2} & 0 & g_1
\end{psmallmatrix}\in \mathrm{GL}(5,\mathbb{C})\;:\; g_1,g_2\in \mathbb{C}, \;g_1\neq 0\right\}.$$

To proceed, we are forced to reduce the computational complexity of the problem. We introduce a mild assumption on the construction of the snake, setting $s_1=s_3$. 

In the end, we verify that the CR structure of the snake with the parameters $s_2=\tfrac12$ and $s_1=s_3$ is described by the vanishing of the four invariants $$\mathbf{J}\equiv 0,\qquad \mathbf{N}\equiv 0,\qquad \mathbf{L}\equiv 0,\qquad \mathbf{F}\equiv 0.$$
The other eight are non-vanishing regardless of the value of $s_1=s_3$.
\bibliographystyle{alpha}
\bibliography{biblio}

\end{document}